\documentclass[a4paper,11pt]{amsart}
\usepackage[colorlinks=true,urlcolor=blue,citecolor=blue,linkcolor=blue,pdfstartview=FitH,pdfpagemode=None,bookmarksnumbered=true]{hyperref}
\usepackage[latin1]{inputenc}
\usepackage[T1]{fontenc}
\usepackage{graphicx}
\usepackage[frenchb]{babel}

\newtheorem{thm}{Th\'{e}or\`{e}me}[section]

\newtheorem{lem}[thm]{Lemme}

\theoremstyle{definition}
\newtheorem{defn}[thm]{D\'{e}finition}

\theoremstyle{remark}
\newtheorem{rem}[thm]{Remarque}
\newtheorem{xple}[thm]{Exemple}

\numberwithin{equation}{section}

\newcommand{\norm}[1]{\left\Vert#1\right\Vert}
\newcommand{\abs}[1]{\left\vert#1\right\vert}
\newcommand{\set}[1]{\left\{#1\right\}}
\newcommand{\Real}{\mathbb R}

\hypersetup{
    pdftitle={Existence de points fixes enlac\'{e}s \`{a} une orbite p\'{e}riodique d'un hom\'{e}omorphisme du plan}, 
    pdfauthor={Christian Bonatti, Boris Kolev}, 
    pdfsubject={MSC : 55M20, 54H20}, 
    pdfkeywords={Hom\'{e}omorphismes du plan, Points fixes}, 
    }

\begin{document}

\title[Existence de points fixes]{Existence de points fixes enlac\'{e}s \`{a} une orbite p\'{e}riodique
d'un hom\'{e}omorphisme du plan}%

\author[C. Bonatti]{Christian Bonatti}
\address{Universit\'{e} de Bourgogne, D\'{e}pt. de Math\'{e}matiques, B.P. 138, 21004 Dijon Cedex, France.}
\email{bonatti@satie.u-bourgogne.fr}%

\author[B. Kolev]{Boris Kolev}%
\address{CMI, 39, rue F. Joliot-Curie, 13453 Marseille cedex 13, France}%
\email{kolev@cmi.univ-mrs.fr}%

\subjclass{55M20, 54H20}%
\keywords{Hom\'{e}omorphismes du plan, Points fixes}%

\date{1 Ao\^{u}t 1991}%


\begin{abstract}
Soit $f$ un hom\'{e}omorphisme du plan qui pr\'{e}serve l'orientation et tel
que $f-Id$ soit une contraction. Sous ces hypoth\`{e}ses, on \'{e}tablit
l'existence, pour toute orbite p\'{e}riodique $\mathcal{O}$, d'un point
fixe ayant un nombre d'enlacement non nul avec $\mathcal{O}$.
\end{abstract}
\maketitle


\section{Introduction}

Un r\'{e}sultat classique d\^{u} \`{a} Brouwer \'{e}nonce que tout hom\'{e}omorphisme du
plan $\Real^{2}$ qui pr\'{e}serve l'orientation et poss\`{e}de une orbite
p\'{e}riodique, poss\`{e}de \'{e}galement un point fixe. Dans le m\^{e}me ordre
d'id\'{e}es, on peut montrer \cite{Gam90,Kol90} qu'un tel hom\'{e}omorphisme
$f$ poss\`{e}de un point fixe li\'{e} \`{a} cette orbite p\'{e}riodique, en ce sens
qu'il n'existe pas de courbe de Jordan $\mathcal{C}$, bordant un
disque $D$ contenant l'orbite p\'{e}riodique mais ne contenant pas le
point fixe, et telle que $f(\mathcal{C})$ soit homotope \`{a}
$\mathcal{C}$ dans le compl\'{e}mentaire de l'orbite p\'{e}riodique et du
point fixe. Une question pos\'{e}e par John Franks dans \cite{Boy89}
demeure toujours sans r\'{e}ponse:

\'{E}tant donn\'{e} un hom\'{e}omorphisme $f:\Real^{2}\rightarrow \Real^{2}$
pr\'{e}servant l'orientation, existe-t-il  pour toute orbite p\'{e}riodique
de $f$ un point fixe ayant un nombre d'enlacement non nul avec cette
orbite p\'{e}riodique?

On sait que la r\'{e}ponse \`{a} cette question est positive pour les
orbites de p\'{e}riode $2$ (voir \cite{Bro90}) ou de p\'{e}riode $3$
\cite{Gua90}. D'autre part, dans une pr\'{e}publication r\'{e}cente, Franks
\cite{Fra92} utilise la r\'{e}ponse affirmative \`{a} cette question comme
\'{e}tant un th\'{e}or\`{e}me de Handel (sans r\'{e}f\'{e}rence): on peut donc supposer
que cette question est, soit r\'{e}solue, soit en passe de l'\^{e}tre.

Nous montrons ici qu'un raisonnement tr\`{e}s simple et tr\`{e}s rapide
permet de r\'{e}pondre par l'affirmative \`{a} la question de Franks pour
les orbites de toutes les p\'{e}riodes des hom\'{e}omorphismes $f$ de
$\Real^{2}$ pour lesquels $f-Id$ v\'{e}rifie un condition de Lipschitz.
Plus pr\'{e}cis\'{e}ment:

\begin{thm}
Soit $f:\Real^{2}\rightarrow \Real^{2}$ un hom\'{e}omorphisme qui
pr\'{e}serve l'orientation et tel que $f-Id$ soit lipschitzienne de
rapport $k\in[0,1]$. Alors, pour toute orbite p\'{e}riodique
$\mathcal{O}=\set{x, f(x), \dotsc , f^{n-1}(x)}$ de $f$, il existe
un point fixe de $f$, $x_{0}$, ayant  un nombre d'enlacement non nul
avec $\mathcal{O}$.
\end{thm}

Outre l'int\'{e}r\^{e}t du r\'{e}sultat, la simplicit\'{e} de la preuve met en
valeur l'avantage qu'il y a \`{a} tester sur cette classe (pas trop
petite) d'hom\'{e}omorphismes, les conjectures concernant les
hom\'{e}omorphismes des surfaces.

\section{Nombre d'enlacement d'un point fixe et d'une orbite p\'{e}riodique}

Soit $f:\Real^{2}\rightarrow \Real^{2}$ un hom\'{e}omorphisme. On note
$Fix(f)$ l'ensemble de ses points fixes. Un point $x\in \Real^{2}$
est p\'{e}riodique de p\'{e}riode $n$ si $f^{n}(x)=x$ mais $f^{k}(x)\neq x$
pour $k\in\set{1,2, \dotsc , n-1}$. On note
\begin{equation*}
    \mathcal{O}(n,f)=\set{x,f(x), \dotsc ,f^{n-1}(x)}
\end{equation*}
l'orbite de $x$ sous $f$.

Soit $\mathrm{x}_{0}\in Fix(f)$, $x$ un point p\'{e}riodique de p\'{e}riode
$n$ et $c$ un arc joignant $x$ et $f(x)$ dans $\Real^{2}\setminus
Fix(f)$, on note $\gamma_{c}$ la courbe ferm\'{e}e obtenue en joignant
bout \`{a} bout les arcs, $c,f(c), \dotsc ,f^{n-1}(c)$.

On note $\omega(x_{0}, \gamma_{c})$ le nombre d'enroulement de
$\gamma_{c}$ autour de $x_{0}$ c'est \`{a} dire le nombre d'intersection
alg\'{e}brique d'une demi-droite g\'{e}n\'{e}rique issue de $x_{0}$ avec
$\gamma_{c}$ (ce nombre est souvent appel\'{e} indice de $x_{0}$ par
rapport \`{a} $\gamma_{c}$).

\begin{lem}\label{lem1}
Soient $c$ et $c^{\prime}$ deux arcs quelconques joignant $x$ et
$f(x)$ dans $\Real^{2}\setminus\set{x_{0}}$. Si $f$ pr\'{e}serve
l'orientation, alors $\omega(x_{0}, \gamma_{c}) - \omega(x_{0},
\gamma_{c^{\prime}})\in n\mathbb{Z}$.
\end{lem}

\begin{proof}
On a
\begin{equation*}
    \omega(x_{0}, \gamma_{c}^{\prime}) = \omega(x_{0},\gamma_{c}) + \sum_{k=o}^{n-1}
    \omega(x_{0},f^{k}(c^{-1}c^{\prime})).
\end{equation*}
Or si $f$ pr\'{e}serve l'orientation
\begin{equation*}
    \omega(x_{0},f^{k}(c^{-1}c^{\prime})) = \omega(x_{0},c^{-1}c^{\prime}),
\end{equation*}
d'o\`{u} $\omega(x_{0}, \gamma_{c}) - \omega(x_{0},\gamma_{c^{\prime}})
= n\omega(x_{0}, c^{-1}c^{\prime})$.
\end{proof}

On peut montrer \'{e}galement que la valeur de
$\omega(x_{0},\gamma_{c})$ (mod $n$) ne d\'{e}pend pas du choix du point
$x$ de $\mathcal{O}$ choisi pour le construire.

\begin{defn}
Avec les notations ci-dessus, on note $Lk(x_{0},\mathcal{O})$
l'unique entier $l\in \set{0,1, \dotsc , n-1}$ tel que
$\omega(x_{0}, \gamma_{c}) -l \in n\mathbb{Z}$ pour un choix
quelconque de $c$ et on l'appelle le \emph{nombre d'enlacement} (ou
\textit{linking number}) du point fixe $x_{0}$ avec l'orbite p\'{e}riodique
$\mathcal{O}$.
\end{defn}

\begin{rem}
L'appellation \og nombre d'enlacement \fg{} est justifi\'{e}e par la remarque
suivante: $Lk(x_{0}, \mathcal{O})$ est aussi le nombre d'enlacement
des deux orbites ferm\'{e}es $C_{x_{0}}$ et $C_{x}$ du champ de vecteurs
canonique induit dans la suspension $T_{f}$ de $f$.
\end{rem}

\begin{xple}
Le point fixe de la rotation d'angle $2k\pi/nk\in\{0, \dotsc ,
n-1\}$ a pour nombre d'enlacement $k$ avec l'une quelconque de ses
orbites p\'{e}riodiques.
\end{xple}

\section{Quelques propri\'{e}t\'{e}s \'{e}l\'{e}mentaires d'un hom\'{e}omorphisme $f$ du plan, tel que $f-Id$ v\'{e}rifie une condition de Lipschitz}

Soit $g:\Real^{2}\rightarrow \Real^{2}$ une application continue. Si
l'ensemble
\begin{equation*}
    \set{k\in[0, \infty[;\forall x, y\in \Real^{2}: \norm{g(x)-g(y)}\leq k \norm{x-y}}
\end{equation*}
est non vide, on note $Lip(g)$ sa borne inf\'{e}rieure. Sinon, on pose
$Lip(g) = +\infty$.

Soient $x$ et $y$ deux points de $\Real^{2}$, on note $[x, y]$ le
segment de droite joignant $x$ et $y$.

\begin{lem}\label{lem2}
Soit $f:\Real^{2}\rightarrow \Real^{2}$ un hom\'{e}omorphisme tel que
$Lip(f-Id)\leq 1$ et $x\in \Real^{2}\setminus Fix(f)$. Alors
$Fix(f)\cap[x,f(x)]=\emptyset$.
\end{lem}

\begin{proof}
Par l'absurde, supposons qu'il existe $y\in Fix(f)\cap[x,f(x)[$. On
a alors:
\begin{equation*}
    \norm{f(x)-x} = \norm{(f(x)-x)-(f(y)-y)} \leq Lip(f-Id) \norm{x-y} <
    \norm{f(x)-x}.
\end{equation*}
D'autre part $f(x)\in \Real^{2}\setminus Fix(f)$ car $f$ est
injective, ce qui conclut.
\end{proof}

\begin{lem}\label{lem3}
Soit $f: \Real^{2}\rightarrow \Real^{2}$ un hom\'{e}omorphisme tel que
$Lip(f-Id) \leq 1$ et $x\in \Real^{2}\setminus Fix(f)$. Alors,
$f([x,f(x)])$ et $[f(x),f^{2}(x)]$ sont homotopes relativement \`{a}
$f(x)$, $f^{2}(x)$ dans $\Real^{2}\setminus Fix(f)$.
\end{lem}

\begin{proof}
Soit $F:[0,1]^{2}\rightarrow \Real^{2}$ d\'{e}finie par:
\begin{equation*}
    F(s, t)=(1-t)((1-s)x+sf(x))+tf((1-s)x+sf(x)).
\end{equation*}
L'image de $F$ est la r\'{e}union des segments de droite
$[y,f(y)](y\in[x,f(x)])$. D'apr\`{e}s le Lemme~\ref{lem2} $y\in[x,f(x)]$
n'est pas un point fixe de $f$, donc $[y,f(y)]\subset
\Real^{2}\setminus Fix(f)$. D'o\`{u} $i\mathrm{m}F\subset
\Real^{2}\setminus Fix(f)$. D'autre part, $F(s,0) = F(0,s) = (1-s)x
+ sf(x)$. On obtient donc une application du disque:
\begin{equation*}
    D=[0,1]^{2}/(s, 0)\sim(0, s)-\Real^{2}\setminus Fix(f)
\end{equation*}
telle que $F(\partial D)=[f(x),f^{2}(x)]\cup f([x,f(x)])$. On a
ainsi r\'{e}alis\'{e} l'homotopie.
\end{proof}

Soient $a$, $b$ deux vecteurs de $\Real^{2}$. On note $(a,
b)\in[0,\pi]$ l'angle non orient\'{e} des vecteurs $a$ et $b$.

\begin{lem}\label{lem4}
Soit $f:\Real^{2}\rightarrow \Real^{2}$ un hom\'{e}omorphisme tel que
$Lip(f-Id)\leq 1$ et $x\in \Real^{2}\setminus Fix(f)$. Alors pour
tout $y\in[x,f(x)]$:
\begin{equation*}
    (f(x)-x,f(y)-y) < \pi/2 .
\end{equation*}
\end{lem}

\begin{proof}
On a:
\begin{equation*}
    \norm{(f(x)-x)-(f(y)-y)} \leq Lip(f-Id)\norm{x-y} \leq
    \norm{f(x)-x},
\end{equation*}
et par suite: $\cos(f(x)-x,f(y)-y)\geq 0$. Si l'angle est $\pi/2$,
cela implique $y=f(y)$ ce qui est impossible en vertu du
Lemme~\ref{lem2}.
\end{proof}

\section{D\'{e}monstration du th\'{e}or\`{e}me}

\begin{figure}
\begin{center}
\includegraphics{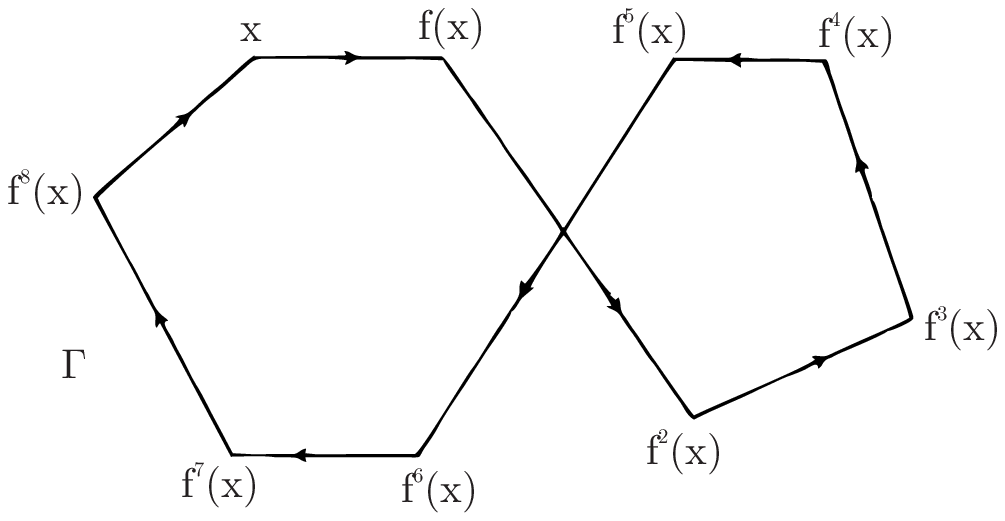}
\caption{} \label{fig1}
\end{center}
\end{figure}

Soit $\mathcal{O} = \set{x,f(x), \dotsc ,f^{n-1}(x)}$ l'orbite d'un
point p\'{e}riodique de p\'{e}riode $n$ d'un hom\'{e}omorphisme de $\Real^{2}$
qui pr\'{e}serve l'orientation et tel que $Lip(f-Id)\leq 1$. Soit
$\Gamma$ la courbe polygonale obtenue en joignant bout \`{a} bout les
segments $[x,f(x)]$, $[f(x),f^{2}(x)]$, $\dotsc$, $[f^{n-1}(x), x]$
(voir Figure~\ref{fig1}) et soit $c=[x,f(x)]$.

\begin{lem}\label{lem5}
Avec les notations ci-dessus, on a:
\begin{enumerate}
  \item $\Gamma\subset \Real^{2}\setminus Fix(f)$,
  \item Pour tout $\mathrm{x}_{0}\in Fix(f)$, $\omega(x_{0}, \gamma_{c}) = \omega(x_{0}, \Gamma)\in\set{-(n-1), \dotsc ,n-1}$.
\end{enumerate}
\end{lem}

\begin{proof}
L'assertion (1) r\'{e}sulte du Lemme~\ref{lem2}. En raisonnant par
r\'{e}currence et en utilisant le Lemme~\ref{lem3}, on \'{e}tablit que
$\Gamma$ est homotope \`{a} $\gamma_{c}$ dans $\Real^{2}\setminus
Fix(f)$, d'o\`{u} l'\'{e}galit\'{e} $\omega(x_{0}, \gamma_{c}) =
\omega(\mathrm{x}_{0}, \Gamma)$. Par ailleurs, le nombre
d'enroulement de $\Gamma$ par rapport \`{a} $x_{0}$ est aussi le nombre
alg\'{e}brique de croisements d'une demi-droite issue de $x_{0}$ avec
$\Gamma$. Ce nombre est donc n\'{e}cessairement inf\'{e}rieur \`{a} $n-1$ en
valeur absolue.
\end{proof}

\begin{figure}
\begin{center}
\includegraphics{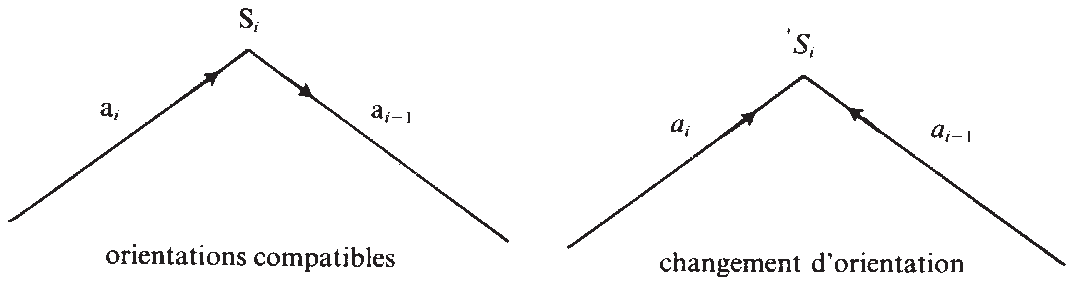}
\caption{} \label{fig2}
\end{center}
\end{figure}

Soient $C_{1},C_{2},\dotsc ,C_{r}$ les composantes connexes born\'{e}es
de $\Real^{2}\setminus\Gamma$ (remarquer qu'il en existe au moins
une, sinon $\Gamma$ serait r\'{e}duit \`{a} un segment de droite ce qui est
exclu en vertu du Lemme~\ref{lem4}) et $C_{\infty}$ la composante
connexe non born\'{e}e. Si $C\in\set{C_{1}, C_{2}\dotsc s, C_{r}}$ on
d\'{e}finit l'indice de $C$ en posant:
\begin{equation*}
    Ind(f, C) = \frac{1}{2\pi} \int_{\partial C}d\varphi ,
\end{equation*}
o\`{u} $\varphi$ d\'{e}signe une d\'{e}termination continue de l'angle du
vecteur $f(x)-x$ avec une direction fixe et o\`{u} $\partial C$ est le
bord orient\'{e} de $C$. Dans chaque composante d'indice non nul il
existe au moins un point fixe de $f$. Soient $S_{0},\dotsc ,
S_{m-1}$ les sommets de $\partial C$ et $a_{0}, \dotsc , a_{m-1}$
ses arr\^{e}tes munies de l'orientation induite par celle de $\Gamma$.
En un sommet $S_{i}$, il y a deux configurations possibles quant \`{a}
l'orientation des ar\^{e}tes adjacentes \`{a} $S_{i}$ : ou bien ces deux
arr\^{e}tes ont des orientations compatibles, ou bien il y a un
changement d'orientation (voir Figure~\ref{fig2}).

Le nombre total de changement d'orientation sur $\partial C$ est un
nombre pair que l'on notera $2p$.

\begin{lem}\label{lem6}
Avec les notations ci-dessus, on a: $Ind(f, C) = 1-p$.
\end{lem}

\begin{figure}
\begin{center}
\includegraphics{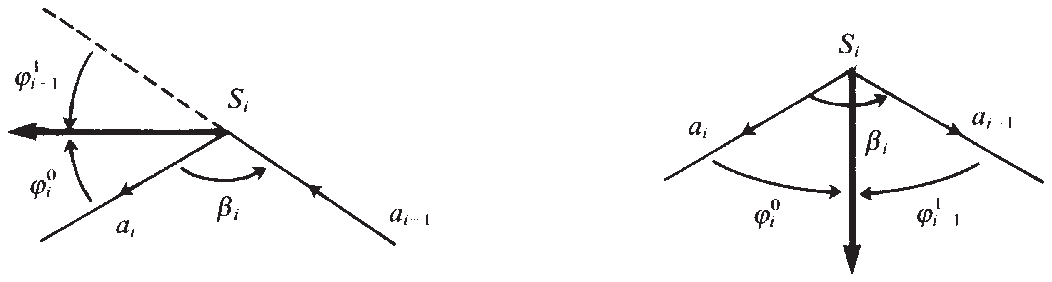}
\caption{} \label{fig3}
\end{center}
\end{figure}

\begin{proof}
On a:
\begin{equation*}
    \int_{\partial C}d\varphi = \sum_{i=0}^{m-1}\int_{a_{i}}d\varphi =
    \sum_{i=0}^{m-1}(\varphi_{i}^{1}-\varphi_{i}^{0}),
\end{equation*}
o\`{u} $\varphi_{i}^{0}$ et $\varphi_{i}^{1}$ sont les valeurs
respectives en $S_{i}$ et $S_{i+1}$ ($i\in \mathbb{Z}/m\mathbb{Z}$)
d'une d\'{e}termination continue $\varphi_{i}$ de l'angle $f(y)-y$
($y\in a_{i}$) avec la tangente orient\'{e}e \`{a} $a_{i}$. Il r\'{e}sulte alors
du Lemme~\ref{lem4} que $\varphi_{i}^{1}-\varphi_{i}^{0}\in[-\pi,
\pi]$ (autrement dit le vecteur $f(y)-y$ ne d\'{e}crit pas de tour
complet lorsque $y$ parcourt $a_{i}$). En d\'{e}signant alors par
$\beta_{i}\in[0, \pi]$ l'angle int\'{e}rieur \`{a} $C$ en $S_{i}$, on a
(voir Figure~\ref{fig3}):
\begin{align*}
    \varphi_{i-1}^{1} - \varphi_{i}^{0} & = \pi-\beta_{i}, \quad
    \text{s'il
n y a pas de changement d'orientation en } S_{i}, \\
    \varphi_{i-1}^{1} - \varphi_{i}^{0} & = -\beta_{i}, \quad
    \text{s'il y a
un changement d'orientation en } S_{i} .
\end{align*}
Par suite:
\begin{align*}
    Ind(f, C) & =
    \frac{1}{2\pi}\left[\sum_{i=0}^{m-1}(\varphi_{i}^{1} - \varphi_{i}^{0})\right]\\
    & =
    \frac{1}{2\pi}\left[\sum_{i=0}^{m-1}(\varphi_{i-1}^{1}-\varphi_{i}^{0})\right]\\
    & = \frac{1}{2\pi}\left[\sum_{i=0}^{m-1}(\pi-\beta_{i})-2p\pi\right]\\
    & = 1-p.
\end{align*}
\end{proof}

Nous pouvons maintenant \'{e}tablir le th\'{e}or\`{e}me. Le nombre d'enroulement
$\omega(x, \Gamma)$ d'un point $x\in \Real^{2}\setminus\Gamma$ ne
d\'{e}pendant que de la composante $C$ \`{a} laquelle il appartient, on
notera $\omega(C, \Gamma)$ cette valeur commune. Il nous reste donc
\`{a} \'{e}tablir l'existence d'une composante $C_{i}$ de
$\Real^{2}\setminus\Gamma$ telle que:
\begin{equation*}
    (1) \quad Ind(f, C_{i})>0 \qquad \text{et} \qquad (2) \quad \omega(C_{i}, \Gamma)\neq 0.
\end{equation*}

\begin{lem}\label{lem7}
Soit $p\in\set{1, \dotsc , r}$ tel que: $\abs{\omega(C_{p},\Gamma)}
= \sup_{i}\abs{\omega(C_{i},\Gamma)}$. Alors $Ind(f, C_{i})>0$.
\end{lem}

\begin{figure}
\begin{center}
\includegraphics{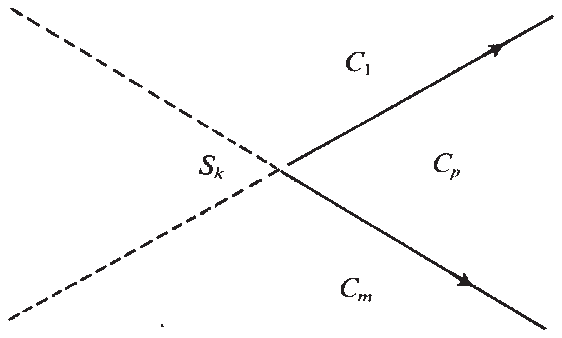}
\caption{} \label{fig4}
\end{center}
\end{figure}

\begin{proof}
Remarquons d'abord que $\omega(C_{\infty}, \Gamma)=0$ et que
$\omega(C_{i}, \Gamma)\neq 0$ pour toute composante $C_{i}$
adjacente \`{a} $C_{\infty}$. Il existe donc bien $p\in\set{1, \dotsc ,
r}$ tel que:
\begin{equation*}
    \abs{\omega(C_{p}, \Gamma)} = \sup_{i}\abs{\omega(C_{i},
\Gamma)}>0.
\end{equation*}
Par l'absurde, supposons que $Ind(f, C_{p}) \leq 0$. D'apr\`{e}s le
Lemme~\ref{lem6}, il existe alors au moins un changement
d'orientation en un des sommets $S_{k}$ de $\partial C_{p}$.

Alors l'une des composantes $C_{l}$ adjacente \`{a} $C_{p}$ en $S_{k}$
v\'{e}rifie $\omega(C_{l}, \Gamma)<\omega(C_{p}, \Gamma)$ et l'autre
$C_{m}$ v\'{e}rifie $\omega(C_{m}, \Gamma)>\omega(C_{p}, \Gamma)$ (voir
Figure~\ref{fig4}). En effet, les valeurs de $\omega(C_{l}, \Gamma)$
et $\omega(C_{m}, \Gamma)$ ne d\'{e}pendent que de l'orientation avec
laquelle on franchit $\Gamma$ pour passer de $C_{p}$ \`{a} $C_{l}$ et
$C_{m}$. Il existe donc $q\in\set{l,m}$ ($q\neq p, \infty$) tel que:
\begin{equation*}
    \abs{\omega(C_{q}, \Gamma)} > \abs{\omega(C_{p}, \Gamma)} >0
\end{equation*}
ce qui contredit l'hypoth\`{e}se faite sur $p$.
\end{proof}

\end{document}